\documentclass[preprint]{imsart}
\RequirePackage[OT1]{fontenc}
\RequirePackage{amsthm,amsmath,natbib}
\RequirePackage[colorlinks]{hyperref}
\RequirePackage{hypernat}
\startlocaldefs
\numberwithin{equation}{section}
\theoremstyle{plain}

\endlocaldefs
\usepackage{psfrag, epsfig, graphicx}
\usepackage{amsfonts,amsmath,latexsym,amssymb}
\usepackage[active]{srcltx}
\usepackage[a4paper]{geometry}
\newtheorem{lemma}{Lemma}
\newtheorem{proposition}{Proposition}

\newtheorem{assumption}{Assumption}
\newtheorem{remark}{Remark}
\newtheorem{definition}{Definition}

\newtheorem{assumption*}{Assumption}
\newtheorem*{theorem}{Theorem}
\renewcommand{\kappa}{\varkappa}

\newcommand{\rd}{{\rm d}}
\newcommand{\cI}{{\cal I}}

\newcommand{\cK}{{\cal K}}
\newcommand{\cL}{{\cal L}}

\newcommand{\cR}{{\cal R}}

\newcommand{\Be}{\boldsymbol{e}}

\newcommand{\Br}{\boldsymbol{r}}

\newcommand{\blb}{\boldsymbol{\beta}}

\newcommand{\bB}{\mathbb B}

\newcommand{\bE}{\mathbb E}
\newcommand{\rE}{\mathbb E}
\newcommand{\bF}{\mathbb F}

\newcommand{\bL}{{\mathbb L}}

\newcommand{\bN}{{\mathbb N}}

\newcommand{\bP}{{\mathbb P}}

\newcommand{\bR}{{\mathbb R}}

\newcommand{\mF}{\mathfrak{F}}

\newcommand{\mL}{\mathfrak{L}}

\newcommand{\mM}{\mathfrak{M}}

\newcommand{\mh}{\mathfrak{h}}

\newcommand{\mz}{\mathfrak{z}}

\newcommand{\epr}{\hfill\hbox{\hskip 4pt
                \vrule width 5pt height 6pt depth 1.5pt}\vspace{0.5cm}\par}

\begin{document}
\begin{frontmatter}
\title{Minimax  estimation of norms of a probability density: II. Rate-optimal estimation procedures}
\runtitle{Estimation of norms of a probability density}
\begin{aug}
\author
{
\fnms{A.} \snm{Goldenshluger}
%\thanksref{m1}
\thanksref{t1}
\ead[label=e1]{goldensh@stat.haifa.ac.il}
}
\and
%\thankstext{t1}{Supported by the ISF grant ???}\ \ \
\author
{\fnms{O. V.} \snm{Lepski}
%\thanksref{m2}
\thanksref{t2}
\ead[label=e2]{oleg.lepski@univ-amu.fr}}
%\thankstext{t1}{Supported by the ISF grant No. 104/11.}\ \ \
%\thankstext{t2}{This work has been carried out in the framework of the Labex Archim\`ede (ANR-11-LABX-0033) and of the A*MIDEX project (ANR-11-IDEX-0001-02), funded by the "Investissements d'Avenir" French Government program managed by the French National Research Agency (ANR).}
%\and
%\author{\fnms{Third} \snm{Author}\thanksref{t1}
%\ead[label=e3]{third@somewhere.com}
%\ead[label=u1,url]{http://www.foo.com}}
%\thankstext{t1}{Some comment}
%\thankstext{t2}{Supported by the ISF grant 389/07}
%\thankstext{t1}{Supported by the ANR grant 0000}
%\thankstext{t3}{Second supporter of the project}
\runauthor{A.~Goldenshluger and O. V.~Lepski}

\affiliation{University of Haifa\thanksmark{m1} and Aix--Marseille Universit\'e, CNRS, Centrale Marseille,
I2M\thanksmark{m2}}

\address{Department of Statistics\\
University of Haifa
\\
Mount Carmel
\\ Haifa 31905, Israel\\
\printead{e1}}

%\affiliation{Aix--Marseille Universit\'e
%Aix Marseille Univ, CNRS, Centrale Marseille, I2M, Marseille, France}

\address{Institut de Math\'ematique de Marseille\\
Aix-Marseille  Universit\'e   \\
 39, rue F. Joliot-Curie \\
13453 Marseille, France\\
%\printead{e1}\\
\printead{e2}}
\end{aug}

\thankstext{t1}{Supported by the ISF grant No. 361/15.}
\thankstext{t2}{This work has been carried out in the framework of the Labex Archim\`ede (ANR-11-LABX-0033) and of the A*MIDEX project (ANR-11-IDEX-0001-02), funded by the "Investissements d'Avenir" French Government program managed by the French National Research Agency (ANR).}

\iffalse
\begin{aug}
\author[t1]{\fnms{A.V.} \snm{Goldenshluger}
\ead[label=e1]{goldensh@stat.haifa.ac.il}}
\thankstext{t1}{Supported by the ISF grant No. 104/11.}\ \ \
\author[t2]
{\fnms{O.V.} \snm{Lepski}
\ead[label=e2]{oleg.lepski@univ-amu.fr}}
%\thankstext{t1}{Supported by the ISF grant No. 104/11.}\ \ \
\thankstext{t2}{This work has been carried out in the framework of the Labex Archim\`ede (ANR-11-LABX-0033) and of the A*MIDEX project (ANR-11-IDEX-0001-02), funded by the "Investissements d'Avenir" French Government program managed by the French National Research Agency (ANR).}
%\and
%\author{\fnms{Third} \snm{Author}\thanksref{t1}
%\ead[label=e3]{third@somewhere.com}
%\ead[label=u1,url]{http://www.foo.com}}
%\thankstext{t1}{Some comment}
%\thankstext{t2}{Supported by the ISF grant 389/07}
%\thankstext{t1}{Supported by the ANR grant 0000}
%\thankstext{t3}{Second supporter of the project}
\runauthor{A.V. Goldenshluger and O.V. Lepski}

\affiliation{University of Haifa and   Aix--Marseille Universit\'e}

\address{Department of Statistics\\
University of Haifa
\\
Mount Carmel
\\ Haifa 31905, Israel\\
\printead{e1}}

\affiliation{Aix Marseille Univ, CNRS, Centrale Marseille, I2M, Marseille, France}

\address{Institut de Math\'ematique de Marseille\\
Aix-Marseille  Universit\'e   \\
 39, rue F. Joliot-Curie \\
13453 Marseille, France\\
%\printead{e1}\\
\printead{e2}}
\end{aug}
\fi

%\maketitle
\begin{abstract}
In this paper we develop rate--optimal estimation procedures in the problem of estimating
the $\bL_p$--norm, $p\in (0, \infty)$ of a probability density  from independent observations.
%This is the second part of the project initiated in \cite{gl20} and related to estimation of $\bL_p$-%norm, $p\in (1,\infty)$, of the probability
%density defined on $\bR^d, d\geq 1$, from independent observations.
The density is assumed to be defined on $\bR^d$, $d\geq 1$ and to belong to a ball in the anisotropic
Nikolskii space.
We adopt the minimax approach and  construct rate--optimal  estimators in the case
of integer $p\geq 2$.
%It is extremely surprising
We demonstrate that,
depending on parameters of Nikolskii's class and the norm index $p$,
the risk asymptotics ranges from inconsistency to $\sqrt{n}$--estimation.
%
%under some special relationships between
%parameters of the Nikolskii class and the norm index $p$ the corresponding
%nonlinear functional admits $\sqrt{n}$-uniformly consistent estimator.
The results in this paper complement the minimax
lower bounds derived in the companion paper \cite{gl20}.
\end{abstract}


\end{frontmatter}

\section{Introduction}
 Suppose that we observe i.i.d. random vectors $X_i\in\bR^d$, $i=1,\ldots, n,$
 with  common probability density $f$.
 Let $p>1$ be a given number. We are interested in estimating
the $\bL_p$-norm of $f$,
 $$
 \|f\|_p:=\bigg[\int_{\bR^d}|f(x)|^p\rd x\bigg]^{1/p},
 $$
from observation $X^{(n)}=(X_1,\ldots,X_n)$.
By estimator of $\|f\|_p$ we mean any $X^{(n)}$-measurable map $\widetilde{N}:\bR^n\to \bR$.
Accuracy of an estimator $\widetilde{N}$
is measured by the quadratic risk
\vskip0.2cm
\centerline{$
 \cR_n[\widetilde{F}, f]:=\Big(\bE_f \big[\widetilde{N}-\|f\|_p\big]^2\Big)^{1/2},
$}
\vskip0.2cm
\noindent where $\bE_f$ denotes  expectation with respect to the probability measure
$\bP_f$ of the observations $X^{(n)}=(X_1,\ldots,X_n)$.
\par
We adopt minimax approach to measuring estimation accuracy.
Let $\mF$ denote the set of all probability densities defined on $\bR^d$.
The maximal risk of an estimator $\widetilde{N}$ on the set
$\bF\subset\mF$ is defined by
%and each estimator $\widetilde{N}$ we associate its the maximal risk on $\bF$ that is
$$
\cR_n\big[\widetilde{N}, \bF\big]:=\sup_{f\in\bF}\cR_n[\widetilde{N}, f],
$$
and the minimax risk is
\[
 \cR_n[\bF]:=\inf_{\widetilde{N}} \cR_n[\widetilde{N}; \bF],
\]
where $\inf$ is taken over all possible estimators.
\par
In the companion paper \cite{gl20} (referred to hereafter as Part~I), we derived lower bounds
on the minimax risk over  functional classes
 $\bF=\bN_{\vec{r},d}\big(\vec{\beta},\vec{L}\big)\cap \bB_q(Q)$,
 where $\bN_{\vec{r},d}\big(\vec{\beta},\vec{L}\big)$ denotes
 anisotropic Nikolskii's class [see Definition~\ref{def:nikolskii} below], and
 $\bB_q(Q):=\{f: \|f\|_q\leq Q\}$
 is the ball in $\bL_q\big(\bR^d\big)$ of radius $Q$.
 %that is
%$$
%\cR_n\big[\bN_{\vec{r},d}\big(\vec{\beta},\vec{L}\big)\cap \bB_q(Q)\big]:=\inf_{\widetilde{N}}\cR_n%\big[\widetilde{N}, \bN_{\vec{r},d}\big(\vec{\beta},\vec{L}\big)\cap \bB_q(Q)\big],
%$$
%where  infimum is taken over all possible estimators.
%Here recall, $\bN_{\vec{r},d}\big(\vec{\beta},\vec{L}\big)$ denotes an anisotropic Nikol'skii class  %(see Definition \ref{def:nikolskii} below)
%and $\bB_q(Q)$ is ball in $\bL_q\big(\bR^d\big)$ of the radius $Q$.
Specifically, we found the sequence $\phi_n$ completely determined by
the $\vec{\beta},\vec{L}, \vec{r}, q, p$ and $n$ such that
\begin{equation*}\label{eq:lower-bound}
\cR_n\big[\bN_{\vec{r},d}\big(\vec{\beta},\vec{L}\big)\cap \bB_q(Q)\big]\gtrsim \phi_n,\; n\to\infty.
\end{equation*}
The goal of the present paper  %present part
is to develop a rate--optimal estimator, say, $\hat{N}$, such that
for any given   $\vec{\beta},\vec{r},\vec{L}, q, Q$
$$
\limsup_{n\to\infty}\,\phi^{-1}_n\,
\cR_n\big[\hat{N}, \bN_{\vec{r},d}\big(\vec{\beta},\vec{L}\big)\cap \bB_q(Q)\big]<\infty.
$$
%Any estimator possessing the above property is called \textsf{minimax} or \textsf{rate-optimal}.
We provide  explicit construction of such an estimator for integer values of
$p\geq 2$.
\par
The problem of estimating nonlinear functionals of a probability density
has been studied in the literature: we refer to Part~I for background and pointers to the relevant
literature. Here we restrict ourselves with brief reminder of   main definitions and
results obtained
in Part~I.
 \par
 We begin with the definition of the anisotropic Nikolskii's classes.
% \paragraph{ Anisotropic Nikol'skii class}
%\label{sec:nikolski}
Let $(\Be_1,\ldots,\Be_d)$ denote the canonical basis of $\bR^d$.
 For a function $G:\bR^d\to \bR^1$ and
real number $u\in \bR$
the first order difference operator with step size $u$ in direction of the variable
$x_j$ is defined  by
$
 \Delta_{u,j}G (x)=G(x+u\mathbf{e}_j)-G(x),\;j=1,\ldots,d.
$
By induction,
the $k$-th order difference operator
%with step size $u$ in direction of the variable $x_j$
is
%defined as
\vskip0.15cm
\centerline{$
 \Delta_{u,j}^kG(x)= \Delta_{u,j} \Delta_{u,j}^{k-1} G(x) = \sum_{l=1}^k (-1)^{l+k}\binom{k}{l}\Delta_{ul,j}G(x).
%g(x+lue_j),\;\;j=1,\ldots,d.
$}
%\vskip-0.15cm
\begin{definition}
\label{def:nikolskii}
For given  vectors $\vec{\beta}=(\beta_1,\ldots,\beta_d)\in (0,\infty)^d$, $\vec{r}=(r_1,$ $\ldots,r_d)\in [1,\infty]^d$,
 and $\vec{L}=(L_1,\ldots, L_d)\in (0,\infty)^d$ a function $G:\bR^d\to \bR^1$ is said to
 belong to  anisotropic
Nikolskii's class $\bN_{\vec{r},d}\big(\vec{\beta},\vec{L}\big)$ if
 $\|G\|_{r_j}\leq L_{j}$ for all $j=1,\ldots,d$ and
 there exist natural numbers  $k_j>\beta_j$ such that
\vskip0.15cm
\centerline{$
 \big\|\Delta_{u,j}^{k_j} G\big\|_{r_j} \leq L_j |u|^{\beta_j},\;\;\;\;
\forall u\in \bR,\;\;\;\forall j=1,\ldots, d.
$}
%\end{itemize}
\end{definition}
Important quantities that are related to Nikolskii's classes and determine
asymptotics of the minimax risk
%and in the estimator construction presented in Sections \ref{sec:ideas} and \ref{sec:main-result}
are the following:
\begin{align*}
&\frac{1}{\beta}:=\sum_{j=1}^d\frac{1}{\beta_j},\quad \frac{1}{\omega}:=\sum_{j=1}^d\frac{1}{\beta_jr_j},\quad L_\beta:=\prod_{j=1}^dL_j^{\frac{1}{\beta_j}},
\\
&\tau(s):=1-\frac{1}{\omega} + \frac{1}{\beta s},\;\;\; s\in[1,\infty].
\end{align*}
 It is worth mentioning that $\tau(\cdot)$ appears in embedding theorems for Nikolskii's spaces;
 see Section~\ref{sec:subsubsec-facts-Nikolskii} for details.
\par
In Part~I we have established the  lower bound (\ref{eq:lower-bound}) on the minimax risk
 with  sequence $\phi_n$ defined by:
\begin{gather}
\theta=\left\{
\begin{array}{clc}
\frac{1}{\tau(1)},\quad&\tau(p)\geq 1;
\\*[2mm]
\frac{1/p-1/q}{1-1/q-(1-1/p)\tau(q)},\quad&\tau(p)< 1,\;\tau(q)<0;
\\*[2mm]
\frac{\tau(p)}{\tau(1)},\quad&\tau(p)<1,\;\tau(q)\geq 0;
\end{array}
\right.
\nonumber\\*[2mm]
\label{eq:rate}
\phi_n=L_\beta^{\frac{1-1/p}{\tau(1)}}n^{-\theta^*},\quad \theta^*=2^{-1}\wedge\theta;
\end{gather}
see Theorem~1 of Part~I. Our goal is to develop a rate--optimal estimator whose risk
converges to zero at the rate $\phi_n$.

 \section{%Ideas led to our
 Estimator construction}
 \label{sec:ideas}
 Assuming that $p$ is integer  let us first discuss the problem of
 estimating  a closely related functional
$\|f\|_p^p$.
\par
Let $K:[-1, 1]^d \to \bR$ be a given function (kernel), and let
 $h=(h_1, \ldots, h_d)\in (0,1]^d$ be a given vector (bandwidth).  Let
$$
K_h(x):=(1/V_h) K(x/h),\quad V_h:=\prod_{k=1}^d h_k,
$$
where here in all what follows $y/x$    denotes   the coordinate--wise division for $x, y\in \bR^d$.
Define
\begin{equation}
\label{eq:bais+smoother}
 S_h(x):= \int K_h(x-y) f(y) \rd y,\;\;\;B_h(x):= S_h(x)-f(x).
\end{equation}
Obviously, $B_h(x)$ is the bias of the kernel density estimator of $f(x)$
associated with kernel~$K$ and bandwidth $h$.
\par
The construction of our estimator for
$
\|f\|_p^p
$
is based on a simple observation formulated below
as Lemma \ref{lem:representation}.
  \begin{lemma}
 \label{lem:representation}
For any  $p\in\bN^*, p\geq 2$,  $f\in\mF$
and  $h\in (0, \infty)^d$ one has
  \begin{eqnarray}
   \|f\|_p^p &=&  (1-p) \int S_h^p(x)\rd x +p\int S_{h}^{p-1}(x) f(x)\rd x
\nonumber
\\
  &&
 + \sum_{j=2}^{p} \tbinom{p}{j} (-1)^j\int  [S_h(x)]^{p-j}B^j_h(x)\rd x.
\label{eq:N-representation}
 \end{eqnarray}
 \end{lemma}
The  proof is elementary and  given in Appendix. It does not require any assumption
on the kernel $K$ except of existence of the integrals on the right hand side of
(\ref{eq:N-representation}).
\par
 Let $\hat{T}_{1,h}$ and $\hat{T}_{2,h}$ be  estimators of
\[
T_p^{(1)}(f):= \int S_h^p(x) \rd x,\;\;\;T_p^{(2)}(f):=\int S_h^{p-1}(x) f(x)\rd x,
\]
respectively.
Then we estimate $\|f\|_p^p$ by
\begin{equation*}%\label{eq:Th}
\hat{T}_h=(1-p) \hat{T}_{1,h} + p\;\hat{T}_{2,h}.
\end{equation*}
Note that if $\hat{T}_{1,h}$ and $\hat{T}_{2,h}$ are unbiased
estimators of $T_p^{(1)}(f)$ and $T_p^{(2)}(f)$ then, in view of Lemma~\ref{lem:representation},
the bias of $\hat{T}_h$ in estimation of $\|f\|_p^p$ is
$$
\sum_{j=2}^{p} \tbinom{p}{j} (-1)^j\int  [S_h(x)]^{p-j}B^j_h(x)\rd x.
$$
The last quantity can be efficiently bounded from above via norms
of the bias  $B_h(\cdot)$ and the underlying density $f$.
\par
The natural unbiased estimators for $T_p^{(1)}(f)$ and $T_p^{(2)}(f)$  are
based on the U--statistics:
\begin{eqnarray*}
 \hat{T}_{1, h} &:=&  \frac{1}{\binom{n}{p}} \sum_{i_1, \ldots, i_p} U^{(1)}_h(X_{i_1}, \ldots, X_{i_p}),
 \\
 \hat{T}_{2, h} &:=&
 \frac{1}{\binom{n}{p}} \sum_{i_1, \ldots, i_p}
 U_h^{(2)}(X_{i_1}, \ldots, X_{i_p}),
\end{eqnarray*}
where the summations are taken
over all possible combinations of $p$ distinct  elements $\{i_1, \ldots, i_p\}$
of $\{1, \ldots, n\}$, and
\begin{eqnarray*}
 U_h^{(1)}(x_1, \ldots, x_p) &:=& \int K_h(y-x_1)\cdots K_h(y-x_p) \rd y,
 \\
 U_h^{(2)}(x_1, \ldots, x_p)&:=& \frac{1}{p}
 \sum_{i=1}^p
 \prod_{\substack{j=1\\ j\ne i}}^{p} K_h(x_{j}-x_{i}).
\end{eqnarray*}
\par
 It is worth mentioning that not only there is the  explicit
 formula for the bias of $\hat{T}_h$, but also its variance admits
a rather simple analytical bound. The following result states an upper bound on the variance of
$\hat{T}_h$; the proof is given in Appendix.
\begin{lemma}\label{lem:var-T1}
Let $K$ be symmetric bounded function supported on $[-1,1]^d$. Then for all $p\in\bN^*, p\geq 2$,  $f\in\mF$ and $h\in (0,\infty)^d$ one has
\begin{align*}
 {\rm var}_f \big[\hat{T}_{h}\big] &\leq  C\|K\|_\infty^{2p} \sum_{k=1}^p
 \|f\|_{2p-k}^{2p-k}\;\big(n^k V_h^{k-1}\big)^{-1},
 \end{align*}
where $C$ is a constant  depending on $p$ and $d$ only.
\end{lemma}
If $\hat{T}_h$ is a "reasonable" estimator of $\|f\|_p^p$ then it seems natural to define
\begin{equation}\label{eq:N-h}
\hat{N}_h:=\big|\hat{T}_h\big|^{1/p}
\end{equation}
 as an estimator for $\|f\|_p$.
%Below we show that in fact $\hat{N}_h$ is rate--optimal in the problem of estimating $\|f\|_p$.
%
\section{Main result}
\label{sec:main-result}
 In this section we demonstrate that the  estimator
$\hat{N}_h$
with properly chosen bandwidth $h$ is a rate--optimal estimator for   $\|f\|_p$,  provided that
 $$
 \vec{r}\in [1,p]^d\cup [p,\infty]^d,\;\;\;p\in \bN^*,\;p\geq 2,
 $$
 i.e., $r_j\leq p$ for all $j=1,\ldots, d$, or $r_j\geq p$ for all $j=1,\ldots, d$.
% whatever $p\in\bN^*$, $p\geq 2$.
The proof is
based on the derivation of  tight uniform upper bounds on the bias and the variance of
 $\hat{T}_h$
over anisotropic Nikolskii's classes.
To get such bounds  for the bias term
we use a special construction of kernel $K$ [see, e.g. \cite{gl14}].
 %\paragraph*{Construction of kernel $K$}
\par
For a given positive integer $\ell\in\bN^*$ and function
 $\cK:\bR\to \bR$ supported on $[-1,1]$  define
\[
\cK_\ell(y)=\sum_{i=1}^\ell \tbinom{\ell}{i} (-1)^{i+1}i^{-1}\cK\big(y/i\big).
\]
\begin{assumption}
\label{ass2:kernel}
$\cK$ is  symmetric, $\int_{\bR} \cK(y)\rd y=1$, $\|\cK\|_\infty<\infty$, and
let
\[
 K(x)=\prod_{j=1}^d \cK_\ell(x_j),\; \forall x\in\bR^d.
\]
\end{assumption}
%
%\paragraph{Construction of minimax estimator}
Let us tntroduce the following notation. For  $j=1,\ldots,d$ let
\begin{equation}\label{eq:p-kappa}
\kappa_j:=\left\{
\begin{array}{ll}
\frac{\beta_j\tau(p)}{\tau(r_j)}, &r_j\leq p,\; \tau(q)>0;
\\*[2mm]
\;\; \beta_j, &\text{otherwise},
\end{array}
\right.
\;\;
p_j:=\left\{
\begin{array}{lll}
\frac{2(1-1/p)}{1-1/r_j}, &r_j\geq p;
\\*[2mm]
\;\;\;2,\; &r_j\leq p, \tau(q)>0;
\\*[2mm]
\frac{2(1/p-1/q)}{1/r_j-1/q}, &r_j<p,\; \tau(q)\leq 0,
\end{array}
\right.
\end{equation}
and let
\[
\frac{1}{\upsilon}:=\sum_{j=1}^d \frac{1}{p_j\kappa_j}.
\]
Define $\mh=(\mh_1,\ldots,\mh_d)$ by
\begin{equation}\label{eq:bandwidth}
\mh_j:=L_j^{-1/\kappa_j}\big(\mL n^{-1}\big)^{\frac{2}{\kappa_j p_j[1+2(1-1/p)/\upsilon]}},
\end{equation}
where  $\mL$ is a constant that is completely determined by the class parameters
$\vec{\beta},\vec{L}, \vec{r}$ and~$p$  (a cumbersome but explicit expression for $\mL$
is given in
Section~ \ref{sec:subsection-proof-theorem}).
%The suggested estimator is then $\hat{N}_\mh$.
\par
The main result of this paper is given in the next theorem.
  \begin{theorem}
 \label{th:p>=3}
Let $p\in\bN^*, p\geq 2$, $Q>0$, $\vec{\beta}\in (0,\infty)^d$, $\vec{L}\in (0,\infty)^d$, $\vec{r}\in [1,p]^d\cup [p,\infty]^d$ and  $q\geq 2p-1$
be fixed.  Let Assumption~\ref{ass2:kernel} hold with $\ell>\max_{j=1,\ldots,d}\beta_j$.
Let $\hat{N}_\mh$ be the estimator (\ref{eq:N-h}) associated with bandwidth
$\mh$ defined in (\ref{eq:bandwidth});
then
$$
\limsup_{n\to\infty}\,\phi^{-1}_n\,\cR_n\big[\hat{N}_\mh, \bN_{\vec{r},d}\big(\vec{\beta},\vec{L}\big)\cap \bB_q(Q)\big]<\infty,
$$
where $\phi_n$ is given in (\ref{eq:rate}).
 \end{theorem}
\begin{remark}\mbox{}
 \par\smallskip
 {\rm (i)} Combining the result of the theorem  with that of Theorem 1 in Part~I
 we conclude that the suggested estimator $\hat{N}_{\mh}$ is rate-optimal.
Thus, the problem of constructing a rate--optimal estimator is solved completely if
 $r_j\leq p$ for all $j=1, \ldots, d$, or $r_j\geq p$ for all $j=1, \ldots, d$.
 %
 %in the dimension 1 or, more generally, in so-called \textsf{semi-isotropic} case, that is
 %$r_j=\boldsymbol{r}, j=1,\ldots, d,$ for some $\boldsymbol{r}\in[1,\infty]$.
%\par
%\noindent
The general setting when values of the coordinates of $\vec{r}$ may be arbitrary
with respect to $p$ remains an open problem.
We conjecture that in the general setting a rate--optimal estimator
does not belong to the family $\{\hat{N}_h, h\in\bR^d\}$,
and a different estimation procedure has to be developed.
Another intriguing question is whether
condition $q\geq 2p-1$ is necessary in order to guarantee the obtained estimation accuracy?
 \par\smallskip
{\rm (ii)}
  It is %extremely
  surprising that
  %such a nonlinear functional as
  the $\bL_p$-norm can be estimated with parametric rate $n^{-1/2}$. To the best of our knowledge
  this phenomenon has not been observed  in the literature. Note that the parametric regime is possible only if $\tau(q)>0$.
   Indeed,
 $$
 \frac{1/p-1/q}{1-1/q-(1-1/p)\tau(q)}\geq \frac{1}{2}
 \;\;\;\Leftrightarrow\;\;
 \frac{2}{p}-1-\frac{1}{q}\geq -\Big(1-\frac{1}{p}\Big)\tau(q),
 $$
 and the last inequality is
impossible if $\tau(q)\leq 0$ because  $p\geq 2$.
\par\smallskip
{\rm (iii)} A particularly simple description of the minimax rate of convergence $\phi_n$
is obtained in the specific case $p=2$ and $q=\infty$.
Here $\theta^*=(\max\{\tau(1),2\})^{-1}$ if $\vec{r}\in [2,\infty]^d$, and
\begin{gather*}
\theta^*=\left\{
\begin{array}{clc}
\tfrac{1}{2},\quad&\omega>1;
\\*[2mm]
\frac{1}{1+1/\omega},\quad&\omega\leq 1,
\end{array}
\right.
\end{gather*}
if $\vec{r}\in [1,2]^d$. As we see, the regime corresponding to the exponent $\tau(p)/\tau(1)$ does not appear if $p=2, q=\infty$.
%Moreover it is easy to see that $\tau(1)\geq 2$ if $\tau(\infty)\leq 0$ and, therefore, independently of whether $\vec{r}\in [2,\infty]^d$ or $\vec{r}\in [1,2]^d$
\end{remark}

\section{Discussion}
In this section we compare and contrast our results with other results in the literature.
Our discussion is restricted
to the problem of estimating the $\bL_p$--norm with integer~$p$ only.
%by an integer $p\geq 2$ because it is the case of this part of our project.
 The fundamental differences between our results and other results
 are mainly due to the following two reasons.
 %
 %We show that the fundamental difference of our results from others is mainly due to the following %reasons.

\begin{itemize}

\item {\em Statistical model.}
%We will try to explain why
We argue that  when estimating the $\bL_p$--norm estimation in the density model,
much better accuracy  can be achieved than in the regression/Gaussian white noise model.
\item {\em Domain of observations.}
We also explain why estimation accuracy is completely different for compactly supported
densities $f$, and densities $f$ supported on the entire space $\bR^d$.
\end{itemize}
\subsection{Norm estimation in density and Gaussian white noise models}
To the best of our knowledge, in the context of the Gaussian white noise model, estimation of
the $\bL_p$--norm was studied only in the one--dimensional case, $d=1$.
So we focus on the one--dimensional Nikolskii class that is
denoted $\bN_{\Br,1}(\blb, \boldsymbol{L})$, $\blb>0,$ $\boldsymbol{L}>0$,
and $\boldsymbol{r}\in[1,\infty]$.
\par
Let $W$ be the standard Wiener process and assume that  we observe the trajectory of
the process
\begin{gather}
\label{eq:white-noise}
X_n(t)=\int_0^t f(u)\rd u +n^{-1/2}W(t), \quad t\in [0,1].
\end{gather}
The estimation of $\|f\|_p$, $1\leq p<\infty$ in the Gaussian white noise
model (\ref{eq:white-noise}) was initiated in \cite{LNS} under assumption that $f$ belongs to H\"older's functional class, i.e.,
$f\in\bN_{\mathbf{\boldsymbol{\infty}},1}(\blb, \boldsymbol{L})$.
The recent paper \cite{{han&co1}} considers the same problem
for  Nikolskii's classes with $\boldsymbol{r}\in [p,\infty)$.
In the Gaussian white noise model (\ref{eq:white-noise})
even in the case of integer $p$
 there is a big difference between   the minimax
rates of convergence  obtained for
odd and even values of $p$.
Such phenomenon does not occur in the density model because density is a positive function.
Thus in the subsequent discussion of the results in the Gaussian white noise model
we restrict ourselves with the case of even $p$.
%speaking about the white gaussian noise model we will always assume that $p$ is an even number.
\par
%Let us come back
According to the results obtained in \cite{LNS} and \cite{{han&co1}} in the case of even $p$,
%The asymptotic of minimax risk found there (expressed in our notations)  is given by
the minimax rate of convergence (expressed in our notation)  is given by
$
\varphi_n=n^{-\frac{1}{1+\tau(1)}}.
$
Note that in both papers $\boldsymbol{r}\geq p$ and, therefore,
$$
\varphi_n=n^{-\frac{1}{1+\tau(1)}}\gg n^{-\frac{1}{\tau(1)}}=\phi_n,
$$
where $\phi_n$ is defined in (\ref{eq:rate}).
Also, since $\tau(1)>1$ for any $\boldsymbol{r}\neq 1$ we conclude that the parametric regime
is impossible in the model (\ref{eq:white-noise}).
Before providing the explanations why the results
discussed above are so different, let us discuss another result.
\par
\cite{LS99} studied a hypothesis testing problem in the model (\ref{eq:white-noise})
when the set of alternatives consists of functions %belonging to
from
$\bN_{\mathbf{\boldsymbol{r}},1}(\blb, \boldsymbol{L})$, $\boldsymbol{r}\in [1,2)$ separated away from zero %function
in the $\bL_2$-norm. It is well known that this problem is equivalent to the problem of estimating  the
$\bL_2$-norm, and the minimax rate of testing coincides %is the same
%as
with the rate of estimation.  The minimax rate found in \cite{LS99} is $\varphi_n=n^{-\frac{1+\tau(2)}{2[1+\tau(1)]}}$ under assumption $\tau(\infty)>0$.
Noting that $\tau(\infty)>0$ implies $q=\infty$ and
 comparing this result with the result obtained in this paper
 for the case $p=2$, $\boldsymbol{r}\in [1,2)$ and $\tau(\infty)>0$ we conclude
$$
\varphi_n=n^{-\frac{1+\tau(2)}{2[1+\tau(1)]}}\gg n^{-1/2}=\phi_n.
$$
The last inequality follows from the fact that $\tau(1)>\tau(2)$.  Thus we
again conclude that the parametric regime is impossible in the model (\ref{eq:white-noise}). Another interesting feature should be mentioned:
the approach used in \cite{LS99} is based on rather sophisticated  pointwise bandwidth selection scheme
while in our estimation procedure the bandwidth is a fixed vector.
\par
The explanation why estimation accuracy in the density model
is better than in the Gaussian white noise model
(\ref{eq:white-noise}) is rather simple.  In fact, the maximal value of the risk in the density model
is attained on the densities with small $\bL_p$-norm.
Analysis of  our estimation strategy shows that in the density model
the variance of
the corresponding $U$-statistics is proportional to the
$\bL_p$-norm of the density
(see Lemma~\ref{lem:var-T1} and Proposition~\ref{prop:variance_p-integer} for details).
%
%that in the density model  the variance of the corresponding $U$-statistics is proportional to the
%$\bL_p$-norm itself.
The smaller this norm, the smaller the stochastic error of the estimation procedure, and
careful bandwidth selection employed in our procedure  takes into account this fact and
improves considerably
the estimation accuracy.
In contrast, in the Gaussian white noise model
(\ref{eq:white-noise})
%this is the model with the additive noise. As the consequence
the stochastic error of the estimation procedure is independent of the signal $f$ and %, therefore,
%remains unchangeable whatever
of the value of its $\bL_p$--norm.

\subsection{Norm estimation
for compactly and non-compactly supported densities}
We start with the following simple observation.
If a probability density $f$ is defined on a compact
set $\cI\subset\bR^d$ then necessarily
$\|f\|_p\geq |\cI|^{-1/p}>0$
for any $p\in (1,\infty]$. This fact reduces estimation of $\|f\|_p$
for compactly supported densities to the problem of
estimating   $\|f\|_p^p$, which is a much smoother functional.
Indeed, let $\widetilde{N}$ be an estimator of $\|f\|_p^p$.
Then, for any density $f$ and any $ p\in\bN^*, p\geq 2$
$$
\Big|\big|\widetilde{N}\big|^{1/p}-\|f\|_p\Big|\leq \|f\|_p^{-(p-1)}\Big|\big|\widetilde{N}\big|-\|f\|^p_p\Big|\leq |\cI|^{1-1/p} \big|\widetilde{N}-\|f\|^p_p\big|.
$$
\par
The problem of estimating $\|f\|_p^p$
has been considered %is very old statistical problem and
by many authors starting from the seminal  paper \cite{bickel-ritov}; see, for instance,
%made fundamental contribution to the minimax and minimax adaptive estimation of this functional, see for instance
\cite{birge-massart},  \cite{picard},   \cite{beatrice}, \cite{cai-low}, \cite{waart} among  others.  It is worth noting that majority of the papers deals with compactly supported densities belonging to a semi--isotropic functional class, that is $r_l=\boldsymbol{r}$ for any $l=1,\ldots,d$. Below we  discuss these results and present them in a unified way.
\par
Let $p\in\bN^*, p\geq 2$ be fixed and let $\bN_{\boldsymbol{r},d}\big(\vec{\beta},\vec{L},\cI\big)$ denote a semi-isotropic Nikolskii class of densities supported on $\cI$. Assume that $\boldsymbol{r}\geq p$. Then the minimax rate of convergence in estimating $\|f\|_p^p$ on this functional class,
and, therefore, of $\|f\|_p$ as well, is given by
\[
\varphi_n=n^{-(\frac{4}{4+1/\beta}\wedge \frac{1}{2})}.
\]
Hence, taking into account that $\tau(1)=1-1/(\boldsymbol{r}\beta)+1/\beta\leq 1+1/(2\beta)$ for any $\boldsymbol{r}\geq p\geq 2$ we conclude that
$$
\varphi_n=n^{-(\frac{4}{4+1/\beta}\wedge \frac{1}{2})} \ll n^{-(\frac{1}{1+1/2\beta}\wedge \frac{1}{2})}= \phi_n.
$$
In particular, the parametric regime in estimating
compactly supported densities is possible
if and only if $\beta\geq 1/4$ which  is a less restrictive condition than
$\tau(1)\leq 2$.

 \section{Proofs}

 \subsection{Preliminaries}
 In this section we collect known facts from functional analysis
 as well as some recent results related to anisotropic Nikolskii's spaces.
 These results will be used in the subsequent proofs.
 \subsubsection{Strong maximal operator}
For locally integrable function $f$ on $\bR^d$
{\em the strong Hardy--Littlewood maximal operator of $f$}
is defined by
\[
 \mM[f](x) = \sup \bigg\{\frac{1}{|I|} \int_{I} f(y)\rd y: x\in I\bigg\},
\]
 where the supremum is taken over all rectangles $I$
with edges parallel to the coordinate axes and containing point $x$; here  $|\cdot|$ stands for  the Lebesgue measure. In view of the Lebesgue differentiation theorem
\begin{equation}\label{eq:maximal-new}
f(x)\leq  \mM[f](x)\quad a.e.
\end{equation}
Moreover,  if $f\in \bL_p(\bR^d)$ then
\begin{equation}\label{eq:maximal-f}
\big\| \mM[f]\big\|_p\leq c_0 \|f\|_p,\;\;\;1<p\leq \infty
\end{equation}
with constant $c_0$ depending on $p$ and $d$ only; see, e.g.,  \cite{Guzman}.

\subsubsection{Some useful inequalities}
For the ease of reference we recall some well known inequalities that are routinely used
in the sequel.
These results can be found, e.g.,  in \cite{Folland}.
\begin{itemize}
\item {\em Interpolation inequality.} Let $1\leq s_0<s<s_1\leq\infty$.
  If $f\in \bL_{s_0}\big(\bR^d\big)\cap \bL_{s_1}\big(\bR^d\big)$ then $f\in \bL_s\big(\bR^d\big)$ and
\[
 \big\|f\big\|_s\leq \big(\big\|f\big\|_{s_0}\big)^{\frac{(s_1-s)s_0}{(s_1-s_0)s}} \big(\big\|f\big\|_{s_1}\big)^{\frac{(s-s_0)s_1}{(s_1-s_0)s}}.
\]
%\item {\em Young's inequality.} If  $f\in \bL_1$ and $g\in \bL_p$, $1\leq p\leq \infty$ then
%\[
% \|f\ast g\|_p\leq \|f\|_1\|g\|_p.
%\]
\item {\em Young's inequality (general form).} Let $1\leq p, q, r\leq \infty$ and
$1+1/r=1/p+1/q$. If $f\in \bL_p\big(\bR^d\big)$ and $g\in \bL_q\big(\bR^d\big)$ then $f\ast g\in \bL_r\big(\bR^d\big)$ and
\[
 \|f\ast g\|_r \leq \|f\|_p \|g\|_q.
\]

%\item {\em Generalized H\"older inequality.} Let $g_l\in\bL_{s_l}\big(\bR^d\big), l=1,\ldots m$, $m\geq 2$, and $\sum_{l=1}^m 1/s_l=1$. Then $\prod_{l=1}^m g_l\in\bL_1\big(\bR^d\big)$ and
%$$
%\Big\| \prod_{l=1}^m g_l\Big\|_1\leq \prod_{l=1}^m \|g_l\|_{s_l}.
%$$
\end{itemize}
\subsubsection{Some facts related to Nikolskii's classes}
\label{sec:subsubsec-facts-Nikolskii}
Let $s>1$, and let the functional class $\bN_{\vec{r}, d}(\vec{\beta}, \vec{L})$ be fixed.
Define for any $j=1, \ldots, d$
 \begin{equation*}
% \label{eq:tau-gamma}
%  \tau(p):=1- \sum_{j=1}^d \frac{1}{\beta_j} \Big(\frac{1}{r_j}- \frac{1}{p}\Big),\;\;
  \gamma_j(s):=\left\{
  \begin{array}{ll}
   \beta_j \tau(s)/\tau(r_j), & r_j<s;\\
   \beta_j, & r_j\geq s,
  \end{array}
\right.
 \end{equation*}
and let
\[
s^*:=s\vee\left[\max_{j=1, \ldots, d}r_j\right], \;\;\vec{s}:=\big(r_1\vee s,\ldots, r_d\vee s\big).
\]
%%%%%%%%%%%%%%%%%%%%%%%%%%%%%%%%%%%%%%%%%%%%%%%%%%%%%%%%%%%%%%%%%%%%%%%
\iffalse
The following result generalizes  embedding theorem
for anisotropic Nikolskii's classes
$\bN_{r, d}(\beta, L)$ proved in \cite[Section~6.9]{Nikolskii}.
\begin{lemma}[\cite{Lepski2015}]
\label{lem:embedding}
 If $\tau(s^*)>0$ then
 \[
 \bN_{\vec{r}, d}(\vec{\beta}, \vec{L}) \subseteq \bN_{\vec{s}, d}(\vec{\gamma}(s), c\vec{L}),
 \]
where constant $c$ is independent of $s$ and $\vec{L}$.
\end{lemma}
\fi
%%%%%%%%%%%%%%%%%%%%%%%%%%%%%%%%%%%%%%%%%%%%%%%%%%%%%%%%%%%%%%%%%%%%%%
Recall that the bias  $B_h(\cdot)=B_h(\cdot,f)$ of a kernel density estimator associated with kernel $K$ is defined in (\ref{eq:bais+smoother}).
The following result has been proved in \cite{gl14}.
\begin{lemma}\label{lem:GL} Let Assumption \ref{ass2:kernel} hold with $\ell>\max_{j=1,\ldots, d} \beta_j$. Then $B_{h}(x)$ admits for any $x\in\bR^d$ the representation $B_{h}(x)=\sum_{j=1}^d B_{h}^{(j)}(x)$
with functions $B_{h}^{(j)}$ satisfying the following inequalities. There exist $C_1>0$ and $C_2>0$ independent of $\vec{L}$
such that for any $f\in \bN_{\vec{r}, d}(\vec{\beta}, \vec{L})$ and $h\in(0, \infty)^d$
\begin{equation}\label{eq:bias-1}
 \big\|B_{h}^{(j)}\big\|_{r_j} \leq C_1 L_j h_j^{\beta_j},\;\;\;\forall j=1, \ldots, d.
\end{equation}
 Moreover, for any $s>1$ satisfying
  $\tau(s^*)> 0$ one has
 \begin{equation}\label{eq:bias-2}
  \big\|B_{h}^{(j)}\big\|_{s_j}\leq C_2 L_j h_j^{\gamma_j(s)},
  \;\;\forall j=1,\ldots, d.
 \end{equation}
 Finally,
 for any $r\in [1,\infty]$ and $R>0$ there exists $C_3>0$ such that for any
 $f\in\bB_r(R)$ and any  $h\in (0, \infty)^d$
 \begin{equation}\label{eq:bias-3}
  \big\|B_{h}^{(j)}\big\|_{r}\leq C_3,
  \;\;\forall j=1,\ldots, d.
 \end{equation}
\end{lemma}
\par
The next lemma presents an inequality
between different norms of a function belonging to anisotropic Nikolskii's class.
The proof of this lemma is given   in Appendix.
\begin{lemma}
\label{lem:inf-norm}
Let $1\leq p<\infty$ be fixed.
For any $s\in (1,\infty]$ satisfying $s\geq p\vee\left[\max_{j=1, \ldots, d}r_j\right]$
and $\tau(s)>0$ there exists
constant $C>0$ independent of  $\vec{L}$ such that  for any  $f\in \bN_{\vec{r}, d}(\vec{\beta}, \vec{L})$
 \begin{equation*}%\label{eq:tilde-L}
  \|f\|_s \leq C  \big(L_{\gamma(s)}\big)^{\frac{(1/p-1/s)\tau(s)}{\tau(p)}}\big(\|f\|_p\big)^{\frac{\tau(s)}{\tau(p)}},\;\;\;L_{\gamma(s)}:= \prod_{j=1}^d L_j^{\frac{1}{\gamma_j(s)}}.
 \end{equation*}
\end{lemma}

\subsection{Reduction to the risk of $\hat{T}_h$}\label{sec:reduction}
We will analyze the risk of the proposed estimator
$\hat{N}_h$ using two different upper bounds that
relate the risk of $\hat{N}_h$ to the risk of $\hat{T}_h$ via
elementary
inequalities:
\begin{eqnarray*}
%\label{eq:bound-1}
 \rE_f\big|\hat{N}_h - \|f\|_p\big|^2 &=& \rE_f \big|\, |\hat{T}_h|^{1/p}- \|f\|_p\big|^2\leq
 \rE_f \big| \hat{T}_h - \|f\|_p^p\big|^{2/p}
 \\
 &\leq&
 \Big[\rE_f\big( \hat{T}_h- \|f\|_p^p\big)^2\Big]^{1/p};
 \nonumber
 \\
 \rE_f\big|\hat{N}_h - \|f\|_p\big|^2 &= &\rE_f \bigg| \frac{\hat{N}_h^p -\|f\|_p^p}{\sum_{i=0}^{p-1} \hat{N}_h^{i} \|f\|_p^{p-1+i}}\bigg|^2
  \leq \frac{\rE_f \big(\hat{T}_h - \|f\|_p^p\big)^2}{\|f\|_p^{2p-2}}.
%\label{eq:bound-2}
  \end{eqnarray*}
Thus, for any underlying density $f$ and any $p\in\bN^*, p\geq 2$
\begin{equation}\label{eq:N-bound}
 \rE_f\big|\hat{N}_h - \|f\|_p\big|^2\leq \min\bigg\{ \Big[\rE_f\big( \hat{T}_h- \|f\|_p^p\big)^2\Big]^{1/p},\;
  \frac{\rE_f \big(\hat{T}_h - \|f\|_p^p\big)^2}{\|f\|_p^{2p-2}} \bigg\}.
\end{equation}

\subsection{Proof of Theorem}
The proof is divided in three steps.
First we establish upper bounds
(uniform over $\bN_{\vec{r},d}(\vec{\beta},\vec{L})\cap \bB_q(Q)$)
on the bias and the variance of the estimator $\hat{T}_h$ for any  value of $h$.
Then using  the derived  upper bounds on the risk of $\hat{T}_h$
and the risk reduction argument given in Section~\ref{sec:reduction}, we
complete the proof.
 \subsubsection{Bound for the bias}
The next proposition states the upper bound on the bias of the estimator $\hat{T}_h$.
\begin{proposition}
\label{prop:bais-p-integer}
Let $\vec{r}\in [1,p]^d\cup [p,\infty]^d$.
Then there exists $C>0$ independent of $\vec{L}$ such that for any $f\in \bN_{\vec{r},d}\big(\vec{\beta},\vec{L}\big)\cap \bB_q(Q)$ and any $h\in (0, \infty)^d$
$$
\bE_f\big|\hat{T}_h-\|f\|_p^p\big|\leq C\|f\|^{p-2}_p\sum_{j=1}^d L^{p_j}_j h_j^{p_j\kappa_j},
$$
where $p_j$ and $\kappa_j$ are defined in (\ref{eq:p-kappa}).
\end{proposition}
\paragraph{Proof}
Our first objective is to prove that
\begin{align}
\label{eq1:proof-p=>3}
 \big|\rE_f\big( \hat{T}_h\big) - \|f\|_p^p\big| \leq c_1\|f\|^{p-2}_p\|B_h\|_p^2, \quad\forall
 h\in (0, \infty)^d.
\end{align}
If $p=2$ then, by Lemma \ref{lem:representation},  (\ref{eq1:proof-p=>3}) holds as equality with $c_1=1$. If $p\geq 3$ then
we have in view of Lemma \ref{lem:representation}
\begin{align}
\label{eq2:proof-p=>3}
 \big|\rE_f\big( \hat{T}_h\big) - \|f\|_p^p\big| \leq  c_2\sum_{k=2}^{p} \int |B_h(x)|^k |S_h(x)|^{p-k}\rd x.
\end{align}
Note that for any $h\in(0, \infty)^d$ in view of (\ref{eq:maximal-new})
$$
|B_h(x)|\leq |S_h(x)|+f(x)=\|K\|_1 \mM[f](x)+f(x)\leq c_3 \mM[f](x)\;\; a.e.
$$
It yields together with (\ref{eq2:proof-p=>3})
\begin{align*}
 \big|\rE_f\big( \hat{T}_h\big) - \|f\|_p^p\big| \leq c_4\int |B_h(x)|^2  \big(\mM[f](x)\big)^{p-2}\rd x.
\end{align*}
Hence, applying H\"older's inequality with exponents $p/2$ and $p/(p-2)$  and (\ref{eq:maximal-f}) we obtain
$$
\big|\rE_f\big( \hat{T}_h\big) - \|f\|_p^p\big| \leq c_4\big\|\mM[f]\big\|_p^{p-2}\|B_h\|_p^2\leq c_5 \|f\|^{p-2}_p\|B_h\|_p^2.
$$
This completes the proof of (\ref{eq1:proof-p=>3}). We deduce from (\ref{eq1:proof-p=>3}) and Lemma~\ref{lem:GL}
\begin{align}
\label{eq3:proof-p=>3}
 \big|\rE_f \hat{T}_h - \|f\|_p^p\big| \leq c_6 \|f\|^{p-2}_p
 \sum_{j=1}^d  \big\|B_h^{(j)}\big\|_p^2,\quad \forall h\in (0, \infty)^d.
\end{align}
\par
Consider now separately three cases. Let $\vec{r}\in[p,\infty]^d$. Then, for any $j=1,\ldots,d,$ applying the interpolation inequality with $s_0=1$, $s=p$ and $s_1=r_j$,
and (\ref{eq:bias-1}) and (\ref{eq:bias-3}) of Lemma \ref{lem:GL}
 we obtain
\begin{align}
\label{eq4:proof-p=>3}
 \big\|B_h^{(j)}\big\|_p^2\leq c_7\big\|B_h^{(j)}\big\|_{r_j}^{\frac{2(1-1/p)}{1-1/r_j}}\leq c_8\big(L_jh_j^{\beta_j}\big)^{p_j}=c_8L_j^{p_j}h_j^{\kappa_jp_j},
\end{align}
for any $f\in\bN_{\vec{r},d}\big(\vec{\beta},\vec{L}\big)$.
\par
Let now  $\vec{r}\in[1,p]^d$ and $\tau(q)>0$, and recall that in this case
$\tau(p)>0$ because
$q>p.$ Applying (\ref{eq:bias-2}) of Lemma \ref{lem:GL} with $s=p$ and $s_j=p, j=1,\ldots,d$, we obtain
\begin{align}
\label{eq5:proof-p=>3}
 \big\|B_h^{(j)}\big\|_p^2\leq c_9\big(L_jh_j^{\gamma_j(p)}\big)^{2}=c_9L_j^{p_j}h_j^{\kappa_jp_j},
\end{align}
for any $f\in\bN_{\vec{r},d}\big(\vec{\beta},\vec{L}\big)$.
\par
It remains to consider the case $\vec{r}\in [1,p)^d$ and $\tau(q)\leq 0$. Applying the interpolation inequality with $s_0=r_j, s=p$ and $s_1=q$
we get for  any $j=1,\ldots,d$ and any  $f\in\bN_{\vec{r},d}\big(\vec{\beta},\vec{L}\big)\cap\bB_q(Q)$
$$
\big\|B_h^{(j)}\big\|^2_p\leq \big(\big\|B_h^{(j)}\big\|_{r_j}\big)^{\frac{2(1/p-1/q)}{1/r_j-1/q}}\big(\big\|B_h^{(j)}\big\|_q\big)^{\frac{2(1/r_j-1/p)}{1/r_j-1/q}}
\leq c_{10}\big\|B_h^{(j)}\big\|_{r_j}^{p_j}.
$$
To get the last inequality we have used (\ref{eq:bias-3}) of Lemma \ref{lem:GL}  with $r=q$ and $R=Q$. It yields together with
 (\ref{eq:bias-1}) of Lemma \ref{lem:GL} that
 \begin{align}
\label{eq6:proof-p=>3}
 \big\|B_h^{(j)}\big\|_p^2\leq c_{11}L_j^{p_j}h_j^{\kappa_jp_j},
\end{align}
The required result  follows from (\ref{eq3:proof-p=>3}), (\ref{eq4:proof-p=>3}), (\ref{eq5:proof-p=>3}) and (\ref{eq6:proof-p=>3}).
\epr

\subsubsection{Bound for the variance} Now we derive
the upper bound on the variance of $\hat{T}_h$.
Recall that $q\geq 2p-1$ and define
\begin{eqnarray*}
&&\cL:=\left\{
\begin{array}{ll}
\displaystyle{\max_{k=1,\ldots, p}}\big(L_{\gamma(2p-k)}\big)^{\frac{(1-k/p)\tau(2p-k)}{\tau(p)}}, & \tau(q)>0;
\\*[2mm]
\quad\quad\quad\quad\quad\quad1, &\tau(q)\leq 0,
\end{array}
\right.
\end{eqnarray*}
\begin{proposition}
\label{prop:variance_p-integer}
Let $\vec{\beta},\vec{L},\vec{r}, q>2p-1$ and $Q>0$ be fixed. Then
there exists $C>0$ independent of $\vec{L}$ such that for any $f\in \bN_{\vec{r},d}\big(\vec{\beta},\vec{L}\big)\cap \bB_q(Q)$
$$
{\rm var}_f \big[\hat{T}_{h}\big] \leq  C\cL\left\{
\begin{array}{lll}
\sum_{k=1}^p
 \big(\|f\|_p\big)^{\frac{(2p-k)\tau(2p-k)}{\tau(p)}}\big(n^k V_h^{k-1}\big)^{-1},\; &\tau(q)>0
 \\*[2mm]
 \sum_{k=1}^p
 \big(\|f\|_p\big)^{\frac{(q-2p+k)p}{q-p}}\;\big(n^k V_h^{k-1}\big)^{-1},\; &\tau(q)\leq 0.
 \end{array}
 \right.
$$
\end{proposition}
\paragraph*{Proof} In view of Lemma \ref{lem:var-T1}
\begin{align*}
&& {\rm var}_f \big[\hat{T}_{h}\big] &\leq  c_1\sum_{k=1}^p
 \|f\|_{2p-k}^{2p-k}\;\big(n^k V_h^{k-1}\big)^{-1},
 \end{align*}
Assume first that $\tau(q)>0$; this  also implies that $\tau(r)>0$ for any $r\leq q.$
Applying Lemma \ref{lem:inf-norm} with $s=2p-k$ we get for any $k\in \{1, \ldots, p\}$
\begin{equation*}
 \big\|f\big\|_{2p-k}^{2p-k}  \leq \big(L_{\gamma(2p-k)}\big)^{\frac{(1-k/p)\tau(2p-k)}{\tau(p)}}\big(\|f\|_p\big)^{\frac{(2p-k)\tau(2p-k)}{\tau(p)}}
%\nonumber \\
 \leq \cL\big(\|f\|_p\big)^{\frac{(2p-k)\tau(2p-k)}{\tau(p)}}.
 %= \cL \big(\|f\|_p\big)^{\frac{(2p-1)\tau(2p-1)}{\tau(p)}-\frac{(k-1)\tau(\infty)}{\tau(p)}}
% \\*[1mm]
 %&=&
 %\cL \big(\|f\|_p\big)^{a+(k-1)b}.
\end{equation*}
\par
Now assume that $\tau(q)\leq 0$. Applying the interpolation inequality with $s_0=p$, $s=2p-k$ and $s_1=q$  we get for any $k\in \{1, \ldots, p\}$
\begin{eqnarray*}
\label{eq8:proof-p=>3}
&& \big\|f\big\|_{2p-k}^{2p-k} \leq \big(\|f\|_p\big)^{\frac{(q-2p+k)p}{q-p}}\big(\|f\|_q\big)^{\frac{(p-k)q}{q-p}}\leq c_2\big(\|f\|_p\big)^{\frac{(q-2p+k)p}{q-p}}.
%\nonumber \\
\end{eqnarray*}
This completes the proof.
\epr
\subsubsection{Completion of the theorem proof}
\label{sec:subsection-proof-theorem}
Now we are in a position to complete the proof of the theorem.
\par
Define
$$
%\frac{1}{\upsilon}=\sum_{j=1}^d \frac{1}{p_j\kappa_j},\quad
N:=\big(\mL n^{-1}\big)^{\frac{1}{1+2(1-1/p)/\upsilon}}, \quad \mL:=\cL^{1/p}L_\kappa^{1-1/p},\quad L_{\kappa}:=\prod_{j=1}^d L_j^{1/\kappa_j},
$$
and recall that $\mh=(\mh_1,\ldots,\mh_d)$ is defined as follows.
$$
\mh_j:=L_j^{-1/\kappa_j}\big(\mL n^{-1}\big)^{\frac{2}{\kappa_j p_j[1+2(1-1/p)/\upsilon]}}=L_j^{-1/\kappa_j}N^{\frac{2}{\kappa_jp_j}}.
$$
\par\smallskip
1$^0$. Several remarks are in order. Direct computations show that
\begin{eqnarray}
\label{eq1:proof-p=2}
N^{p-2}\sum_{j=1}^d L^{p_j}_j \mh_j^{p_j\kappa_j}&=&dN^p;
\\
\label{eq2:proof-p=2}
\cL N^{p}\big[n^{p}(V_\mh)^{p-1}\big]^{-1}&=&N^{2p}.
\end{eqnarray}
The following useful equality is deduced from (\ref{eq2:proof-p=2}):
\begin{equation}
\label{eq10:proof-p=>3}
\big[n V_\mh\big]^{-1}=\big[\cL^{-1}nN^{p}\big]^{\frac{1}{p-1}}=\big[\cL^{-1}L_{\kappa}\big]^{1/p}N^{1-\frac{2}{p\upsilon}}.
\end{equation}
\par
Let $R_n(f)$ denote the quadratic risk of the estimator $\hat{T}_\mh$, and let
\begin{eqnarray*}
\bF:= \bN_{\vec{r}, d}(\vec{\beta}, \vec{L})\cap\bB_q(Q)\cap \bB_p(N),\quad \bar{\bF}:=\big[\bN_{\vec{r}, d}(\vec{\beta}, \vec{L})\cap\bB_q(Q)\big]\setminus\bF
\end{eqnarray*}
We obviously have
$$
\cR_n:=\cR_n\big[\hat{N}_\mh, \bN_{\vec{r},d}\big(\vec{\beta},\vec{L}\big)\cap \bB_q(Q)\big]=\max\Big(\sup_{f\in\bF}\cR_n\big[\hat{N}_\mh,f\big],
\sup_{f\in\bar{\bF}}\cR_n\big[\hat{N}_\mh,f\big]\Big)
$$
and,
we deduce from (\ref{eq:N-bound}) that
\begin{equation*}
%\label{eq4:proof-p=2}
\cR^2_n\leq \max\Big(\sup_{f\in\bF}\big[R_n(f)\big]^{1/p},\;
\sup_{f\in\bar{\bF}}\|f\|_p^{2-2p}R_n(f)\Big)
\end{equation*}
\par\smallskip
2$^0$.
In view of  Propositions~\ref{prop:bais-p-integer}
and~\ref{prop:variance_p-integer}, and by (\ref{eq1:proof-p=2}) and (\ref{eq10:proof-p=>3})
we have that for any $f\in\bF$
\begin{gather}
\label{eq11:proof-p=>3}
\bE_f\big|\hat{T}_h-\|f\|_p^p\big|\leq c_1 N^p;
\\*[3mm]
\nonumber
{\rm var}_f \big[\hat{T}_{h}\big] \leq   \frac{C_1\big(\vec{L}\big)}{n}\left\{
\begin{array}{lll}
\sum_{k=1}^p
 N^{\frac{(2p-k)\tau(2p-k)}{\tau(p)}+(1-\frac{2}{p\upsilon})(k-1)},\; &\tau(q)>0
 \\*[2mm]
 \sum_{k=1}^p
 N^{\frac{(q-2p+k)p}{q-p}+(1-\frac{2}{p\upsilon})(k-1)},\; &\tau(q)\leq 0.
 \end{array}
 \right.
\end{gather}
Setting $\theta=1-\frac{2}{p\upsilon}+b$,
$$
a=\left\{
\begin{array}{cc}
\frac{(2p-1)\tau(2p-1)}{\tau(p)}, & \tau(q)>0;
\\*[2mm]
\frac{(q-2p+1)p}{q-p}, &\tau(q)\leq 0,
\end{array}
\right.
\quad\;
b=\left\{
\begin{array}{cc}
-\frac{\tau(\infty)}{\tau(p)}, & \tau(q)>0;
\\*[2mm]
\frac{p}{q-p}, &\tau(q)\leq 0,
\end{array}
\right.
$$
we obtain %\footnote{Here and later $y_-=\min(y,0)$}
for any $f\in\bF$
\begin{equation*}
\label{eq13:proof-p=>3}
{\rm var}_f \big[\hat{T}_{h}\big] \leq \frac{C_2\big(\vec{L}\big)N^a}{n}\;\sum_{m=0}^{p-1}
 N^{(1-\frac{2}{p\upsilon}+b)m}=\frac{C_2\big(\vec{L}\big)N^{a+[\theta(p-1)]_-}}{n}=C_3\big(\vec{L}\big)N^{\mz},
\end{equation*}
where  here and later $y_-:=\min(y,0)$, and
$\mz:=a+1+\frac{2(p-1)}{p\upsilon}+[\theta(p-1)]_-$.
It yields together with (\ref{eq11:proof-p=>3})
\begin{equation}
\label{eq14:proof-p=>3}
\sup_{f\in\bF}\big[R_n(f)\big]^{1/p}\leq C_4\big(\vec{L}\big)\Big[N^2+N^{\mz/p}\Big].
\end{equation}
\par\smallskip
3$^0$. Let us compute the quantity $1/\upsilon$. If $\vec{r}\in[p,\infty]^d$ we have
\begin{equation}
\label{eq141:proof-p=>3}
\frac{1}{\upsilon}=\frac{p}{2(p-1)}\sum_{j=1}^d \frac{1-1/r_j}{\beta_j}=\frac{p(1/\beta-1/\omega)}{2(p-1)}.
\end{equation}
If $\vec{r}\in[1,p]^d$ and $\tau(q)>0$ we have
\begin{equation}
\label{eq142:proof-p=>3}
\frac{1}{\upsilon}=\frac{1}{2\tau(p)}\sum_{j=1}^d\frac{\tau(r_j)}{\beta_j}=\frac{1}{2\tau(p)}\bigg(\frac{\tau(\infty)}{\beta}+\frac{1}{\omega\beta}\bigg)
=\frac{1}{2\beta\tau(p)}.
\end{equation}
Finally, if $\vec{r}\in[1,p]^d$ and $\tau(q)\leq 0$ we have
\begin{equation}
\label{eq143:proof-p=>3}
\frac{1}{\upsilon}=\frac{pq}{2(q-p)}\sum_{j=1}^d\frac{1/r_j-1/q}{\beta_j}=\frac{pq}{2(q-p)}\bigg(\frac{1}{\omega}-\frac{1}{\beta q}\bigg).
%=\frac{q}{q-2}-\frac{q\tau(q)}{q-2}\geq \frac{q}{q-2}.
\end{equation}
It yields, in particular,
$$
(p-1)\theta=\left\{
\begin{array}{ccc}
\frac{p-1}{p\beta\tau(p)}-\frac{1}{\beta}+\frac{1}{\omega}, &\vec{r}\in[p,\infty]^d;
\\*[2mm]
0, &\vec{r}\in[1,p]^d,\;&\tau(q)> 0;
\\*[2mm]
\frac{(p-1)q\tau(q)}{q-p}, &\vec{r}\in[1,p]^d,\;&\tau(q)\leq 0.
\end{array}
\right.
$$
Noting that $\vec{r}\in[p,\infty]^d$ implies $\tau(p)\geq 1$  we assert
$$
\frac{p-1}{p\beta\tau(p)}-\frac{1}{\beta}+\frac{1}{\omega}\leq \frac{p-1}{p\beta}-\frac{1}{\beta}+\frac{1}{\omega}=\frac{1}{\omega}-\frac{1}{p\beta}=1-\tau(p)\leq 0.
$$
Hence in all cases $(p-1)\theta\leq 0$ and we conclude that
\begin{eqnarray*}
\mz&=&a+1+\frac{2(p-1)}{p\upsilon}+\theta(p-1)=a+b(p-1)+p
\\
&=&p+\frac{(2p-1)\tau(\infty)+1/\beta}{\tau(p)}-\frac{(p-1)\tau(\infty)}{\tau(p)}=2p.
\end{eqnarray*}
Hence we deduce from (\ref{eq14:proof-p=>3}) that
\begin{equation}
\label{eq144:proof-p=>3}
\sup_{f\in\bF}\big[R_n(f)\big]^{1/p}\leq C_4\big(\vec{L}\big)N^2.
\end{equation}
\par\smallskip
4$^0$.  Note that for any $f\in\bar{\bF}$ one has in view of Proposition \ref{prop:bais-p-integer} and (\ref{eq1:proof-p=2})
\begin{equation}
\label{eq15:proof-p=>3}
\|f\|_p^{2-2p}\big(\bE_f\big|\hat{T}_h-\|f\|_p^p\big|)\big)^2\leq c_1 \|f\|_p^{-2}\bigg(\sum_{j=1}^d L^{p_j}_j h_j^{p_j\kappa_j}\bigg)^2\leq c_2 N^{2}.
\end{equation}
Here we have used that $N$ is small then $n$ is large.
\par
(a). If $\tau(q)>0$  then for any $k\in\{1,\ldots, p\}$
\begin{eqnarray*}
&&(2p-k)\tau(2p-k)-(2p-2)\tau(p)=-[(k-2)\tau(\infty)+(p-2)/(p\beta)]
\\
&&=-(k-2)\big[\tau(\infty)+1/(p\beta)\big]-(p-k)/(p\beta).
\\
&&=-(k-2)\tau(p)-(p-k)/(p\beta).
\end{eqnarray*}
Since $\tau(p)>0$,  for any $k\in\{2,\ldots, p\}$
$$
\frac{(2p-k)\tau(2p-k)}{\tau(p)}-(2p-2)=-(k-2)-\frac{(p-k)}{p\beta\tau(p)}\leq 0.
$$
Putting $\alpha=1-\frac{(p-1)}{\tau(p)p\beta}$ we deduce from Proposition \ref{prop:variance_p-integer} for any $f\in\bar{\bF}$
\begin{eqnarray*}
{\rm var}_f \big[\hat{T}_{h}\big]&\leq& C\cL n^{-1}\bigg[N^{\alpha_-}+\sum_{k=2}^p
N^{-(k-2)-\frac{(p-k)}{p\beta\tau(p)}}\big(n V_h\big)^{1-k}\bigg]
\nonumber\\
&\leq&C_4\big(\vec{L}\big)n^{-1}\bigg[N^{\alpha_-}+\sum_{k=2}^p
N^{-(k-2)-\frac{(p-k)}{p\beta\tau(p)}+(1-\frac{2}{p\upsilon})(k-1)}\bigg]
\nonumber\\
&\leq&C_5\big(\vec{L}\big)n^{-1}\big[N^{\alpha_-}+
N^{\gamma+(p-2)\theta}\big],
\end{eqnarray*}
where we have put $\gamma=1-\frac{2}{p\upsilon}-\frac{p-2}{p\beta\tau(p)}$ and used that $\theta\leq 0$.
The last bound can be rewritten as
\begin{eqnarray}
\label{eq16:proof-p=>3}
{\rm var}_f \big[\hat{T}_{h}\big]
&\leq&C_6\big(\vec{L}\big)\big[n^{-1}N^{\alpha_-}+
N^{z}\big],
\end{eqnarray}
where $z=\gamma+1+\frac{2(p-1)}{p\upsilon}+(p-2)\theta$. We have
\begin{equation*}
z=p-\frac{p-2}{p\beta\tau(p)}+(p-2)b=p-\frac{p-2}{\tau(p)}\big[\tau(\infty)+1/(p\beta)\big]=2
\end{equation*}
and (\ref{eq16:proof-p=>3}) becomes
\begin{eqnarray}
\label{eq17:proof-p=>3}
{\rm var}_f \big[\hat{T}_{h}\big]
&\leq&C_7\big(\vec{L}\big)\big[n^{-1}N^{\alpha}+
N^{2}\big],\quad \forall f\in\bar{\bF}.
\end{eqnarray}
Thus, if $\tau(q)>0$ we conclude from (\ref{eq144:proof-p=>3}), (\ref{eq15:proof-p=>3}) and (\ref{eq17:proof-p=>3})
\begin{eqnarray}
\label{eq18:proof-p=>3}
\cR^2_n\leq C_7\big(\vec{L}\big)\big[n^{-1}N^{\alpha}+
N^{2}\big].
\end{eqnarray}
Additionally we remark that if $\alpha<0$ then
$$
n^{-1}N^{\alpha_-}=C_8\big(\vec{L}\big)N^{\alpha+1+2(1-1/p)/\upsilon}.
$$
If $\vec{r}\in [p,\infty]^d$ then $\tau(p)>1$, and we get from (\ref{eq141:proof-p=>3})
$$
\alpha+1+2(1-1/p)/\upsilon=2-\frac{(p-1)}{\tau(p)p\beta}+\frac{1}{\beta}-\frac{1}{\omega}\geq 2+\frac{1}{p\beta}-\frac{1}{\omega}\geq 2.
$$
If $\vec{r}\in [1,p]^d$ then in view of (\ref{eq142:proof-p=>3})
$$
\alpha+1+2(1-1/p)/\upsilon=2-\frac{(p-1)}{\tau(p)p\beta}+\frac{(p-1)}{p\beta\tau(p)}=2.
$$
Thus, (\ref{eq18:proof-p=>3}) is equivalent to
\begin{eqnarray}
\label{eq19:proof-p=>3}
\cR^2_n\leq C_9\big(\vec{L}\big)\left\{
\begin{array}{ccc}
n^{-1}\vee N^2, &\alpha\geq 0;
\\*[2mm]
N^2, &\alpha<0.
\end{array}
\right.
\end{eqnarray}
The definition of $N$ implies that
$$
N^2\leq n^{-1}\;\Leftrightarrow\; 2(1-1/p)/\upsilon\leq 1\;\Leftrightarrow\; 1\geq \left\{
\begin{array}{ccc}
1/\beta-1/\omega, &\vec{r}\in [p,\infty]^d;
\\*[2mm]
\frac{(p-1)}{p\beta\tau(p)}, &\vec{r}\in [1,p]^d.
\end{array}
\right.
$$
In the case $\vec{r}\in [1,p]^d$ it yields immediately that $N^2\leq n^{-1}\Leftrightarrow a\geq 0$.

It remains to note that if $\vec{r}\in[p,\infty]^d$ then

\vskip0.2cm

\centerline{$
\alpha\geq 0\Leftrightarrow 1-\frac{(p-1)}{p\beta\tau(p)}=1-\frac{1}{\beta}+\frac{1}{\omega}+\bigg[\frac{(p-1)}{p\beta}-\frac{(p-1)}{p\beta\tau(p)}\bigg]
+\bigg[\frac{1}{p\beta}-\frac{1}{\omega}\bigg]\geq 0.
$}

\vskip0.1cm

\noindent
Since in the considered case $\tau(p)\geq 1\Leftrightarrow 1/(p\beta)>1/\omega$ we assert that
$$
1\geq 1/\beta-1/\omega\; \Rightarrow\; a\geq 0.
$$
Thus, we deduce from (\ref{eq19:proof-p=>3})

\vskip0.2cm

\centerline{$
\cR^2_n\leq C_9\big(\vec{L}\big)\max\Big[n^{-1},n^{-\frac{2}{1+2(1-1/p)/\upsilon}}\Big]
$}

\vskip0.2cm

 \noindent and the assertion of the theorem in the case $\tau(q)>0$ follows from (\ref{eq141:proof-p=>3}) and (\ref{eq142:proof-p=>3}).

\par\smallskip

(b).~ Let $\tau(q)\leq 0$. For any $k\in\{1,\ldots p\}$ one has
\[
\frac{(q-2p+k)p}{q-p}-2(p-1)=2-p+\frac{p(k-p)}{q-p}\leq 0
\]
and, therefore, we have in view of Proposition \ref{prop:variance_p-integer} and (\ref{eq10:proof-p=>3})  for any $f\in\bar{\bF}$
\begin{eqnarray*}
{\rm var}_f \big[\hat{T}_{h}\big] &\leq&  C\cL N^{2-p}
 \sum_{k=1}^p
 N^{\frac{(k-p)p}{q-p}}\;\big(n^k V_h^{k-1}\big)^{-1}
 \\
 &\leq& C\cL N^{\frac{2q-p-pq}{q-p}}n^{-1}
 \sum_{k=1}^p
 \big(N^{\frac{p}{q-p}}(n V_h)^{-1}\big)^{k-1}
 \\
  &=& C_{10}\big(\vec{L}\big) N^{2-\frac{q(p-1)\tau(q)}{q-p}}
 \sum_{k=1}^p
 N^{\frac{q\tau(q)(k-1)}{q-p}}= C_{10}\big(\vec{L}\big)N^2.
\end{eqnarray*}
To get the penultimate equality we have used (\ref{eq143:proof-p=>3}).
It yields together with (\ref{eq15:proof-p=>3})
\[
\displaystyle{\sup_{f\in\bar{\bF}}\|f\|_p^{2-2p}R_n(f)\leq  C_{11}\big(\vec{L}\big)N^2}.
\]
which together with (\ref{eq144:proof-p=>3}) leads to
\[
\cR^2_n\leq  C_{11}\big(\vec{L}\big)N^2=C_{11}\big(\vec{L}\big)n^{-\frac{2}{1+2(1-1/p)/\upsilon}}=C_{11}\big(\vec{L}\big)n^{-\frac{2(1/p-1/q)}{1-1/q-(1-1/p)\tau(q)}}.
\]
This completes the theorem proof.
\epr

\section{Appendix}

\paragraph{Proof of Lemma \ref{lem:representation}}
Using the  identity
$
a^p=\sum_{j=0}^p \tbinom{p}{j}b^{p-j}(a-b)^j
$
with $a=f(x)$ and $b=S_h(x)$ we obtain for all $h\in (0, \infty)^d$ and $x\in\bR^d$
\begin{eqnarray*}
  f^p(x)&= &\sum_{j=0}^p \tbinom{p}{j}  [S_h(x)]^{p-j}  \big[f(x)-S_h(x)\big]^j\rd x
  \nonumber
  \\
  &=& S^p_h(x)(1-p)+p[S_h(x)]^{p-1}f(x)+\sum_{j=2}^p \tbinom{p}{j} (-1)^j [S_h(x)]^{p-j}B^j_h(x).
 \end{eqnarray*}
 Integrating the last equality, we come to the statement  of the lemma.
 \epr

\paragraph{Proof of Lemma \ref{lem:var-T1}}
For any $y\in \bR^d$ let
$$
I_h(y):=\otimes_{j=1}^d\big\{x\in \bR^d: |x_j-y_j|\leq h_j\big\}.
$$
\par\smallskip
1$^0$. We start with bounding the variance of $\hat{T}_{1,h}$.
Define
\begin{align*}
 g_k(x_1, \ldots, x_k) &:= \rE_f \big[U_h^{(1)}(x_1, \ldots, x_k, X_{k+1}, \ldots, X_p)\big],\;\;\;k=1,\ldots, p-1,
\\
 g_p(x_1, \ldots, x_p) &:=U_h^{(1)}(x_1, \ldots, x_p).
 \end{align*}
 Then the variance of $\hat{T}_{1,h}$ is
given by the following well known formula [see, e.g., \cite{Serfling}]:
\[
 {\rm var}_f\big[ \hat{T}_{1,h}\big] = \frac{1}{\binom{n}{p}}
 \sum_{k=1}^p \binom{p}{k} \binom{n-p}{p-k} \zeta_k,\;\;\;
 \zeta_k:={\rm var}_f\big[g_k(X_1, \ldots, X_k)\big].
\]
We note that $g_k(x_1, \ldots, x_k)$, $k=1, \ldots, p$ are symmetric functions.
Observe that for $k=1, \ldots, p-1$
\begin{align*}
 |g_k(x_1, \ldots, x_k)| &= \big|\rE_f\big[ g(x_1, \ldots, x_k, X_{k+1}, \ldots, X_p)\big]\big|
 \\
 &\leq \int \prod_{i=1}^k |K_h(y-x_i)| \Big[\prod_{i=k+1}^p \int |K_h(y-x_i)|f(x_i)\rd x_i\Big] \rd y
 \\
 & \leq \|K\|_\infty^{p-k}
 \int \prod_{i=1}^k |K_h(y-x_i)|  \big(\mM[f](y)\big)^{p-k}\rd y.
\end{align*}
Then
\begin{eqnarray*}
 \zeta_k &\leq& \rE_f \big[g_k(X_1, \ldots, X_k)\big]^2
\\
&\leq& \|K\|_\infty^{2(p-k)} \iint \big(\mM[f](y)\big)^{p-k}\big(\mM[f](z)\big)^{p-k}
\\
&&\;\;\;\times\;
\Big[ \prod_{i=1}^k \int K_h(y-x_i)K_h(z-x_i) f(x_i) \rd x_i \Big]\rd y \rd z
\\
&\leq& \|K\|_\infty^{2p}V_h^{-k} \iint \big(\mM[f](y)\big)^{p-k}\big(\mM[f](z)\big)^{p-k}
{\bf 1}\{y-x\in I_{2h}(0)\}
\\
&&\;\;\;\times\;\Big[ \prod_{i=1}^k\int K_h(z-x_i) f(x_i) \rd x_i \Big] \rd y \rd z
\\
&\leq&  \|K\|_\infty^{2p}V_h^{-k}  \iint \big(\mM[f](y)\big)^{p-k}\big(\mM[f](z)\big)^{p}
{\bf 1}\{y-x\in I_{2h}(0)\}  \rd z\rd y
\\
&\leq& c_1\|K\|_\infty^{2p}V_h^{-k+1} \int  \big(\mM[f](y)\big)^{p-k} \mM\big[\mM^p[f]\big](y)\rd y.
 \end{eqnarray*}
Furthermore, by the H\"older inequality and  (\ref{eq:maximal-f})
\begin{align}
 &\int \big(\mM[f](y)\big)^{p-k} \mM\big[\mM^p[f]\big](y) \leq
 \big\|\mM[f]\big\|_{2p-k}^{p-k} \big\| \mM\big[\mM^p[f]\big]\big\|_{(2p-k)/p}
 \nonumber
 \\*[2mm]
 &\leq c_2 \|f\|_{2p-k}^{p-k} \big\|\mM^p[f]\big\|_{(2p-k)/p}
%\nonumber
% \\*[2mm]
 =
c_2 \|f\|_{2p-k}^{p-k} \big\|\mM[f]\big\|_{2p-k}^p \leq c_3 \big\|f\big\|_{2p-k}^{2p-k}.
 \label{eq:mu-M}
 \end{align}
Therefore we obtain for any $k=1,\ldots, p-1$
\[
 \zeta_k \leq c_2 \|K\|_\infty^{2p}V_h^{-k+1} \int [f(x)]^{2p-k}\rd x.
\]
In addition,
\begin{align*}
 \rE_f [g_p(X_1, \ldots, X_p)]^2 =\iint
 \Big[ \prod_{i=1}^p \int K_h(y-x_i)K_h(z-x_i) f(x_i) \rd x_i \Big]\rd y \rd z
 \\
 \leq \frac{\|K\|_\infty^{2p}}{V_h^p}
 \iint {\bf 1}\{y-z\in I_{2h}(0)\} \big(\mM[f](y)\big)^{p}\rd y \rd z  \leq
 \frac{c_3\|K\|_\infty^{2p}}{V_h^{p-1}} \int [f(x)]^p\rd x.
\end{align*}
Thus we obtain
\begin{align}
 {\rm var}_f \big[\hat{T}_{1,h}\big] &\leq  C_1\|K\|_\infty^{2p} \sum_{k=1}^p
 \frac{1}{n^k V_h^{k-1}} \int [f(x)]^{2p-k}\rd x.
\label{eq:var-T1}
 \end{align}
\par\smallskip
 2$^0$.
Bounding  the variance of $\hat{T}_{2,h}$ goes along the same lines. Define
\begin{align*}
 g_k(x_1, \ldots, x_k) &:= \rE_f \big[U_h^{(2)}(x_1, \ldots, x_k, X_{k+1}, \ldots, X_p)\big],\;\;\;k=1,\ldots, p-1,
\\
 g_p(x_1, \ldots, x_p) &:=U_h^{(2)}(x_1, \ldots, x_p).
 \end{align*}
We have for $k=1, \ldots, p-1$
\begin{align*}
 &g_k(x_1, \ldots, x_k)=\rE_f\big[U^{(2)}_h(x_1, \ldots, x_k, X_{k+1}, \ldots, X_p)\big]
 \\
& = \frac{1}{p}\sum_{i=1}^k \prod_{\substack{j=1\\ j\ne i}}^k K_h(x_j-x_i)
 \bigg[\prod_{j=k+1}^p \int K_h(x_j-x_i) f(x_j)\rd x_j\bigg]&
 \\
&\;\;\; +\frac{1}{p} \sum_{i=k+1}^p
 \int
 \prod_{j=1}^k K_h(x_j-x_i)
% \\
% &&\;\;
 \bigg[\prod_{\substack{j=k+1\\j\ne i}}^p
 \int K_h(x_j-x_i) f(x_j)\rd x_j \bigg] f(x_i)\rd x_i.&
\end{align*}
Therefore
\begin{align*}
& |g_k(x_1, \ldots, x_k)| \leq \frac{\|K\|_\infty^{p-k}}{p}\sum_{i=1}^k
\big(\mM[f](x_i)\big)^{p-k} \prod_{\substack{j=1\\ j\ne i}}^k |K_h(x_j-x_i)|
 \\
 &\;\;\;+ \frac{\|K\|_\infty^{p-k-1}}{p}\sum_{i=k+1}^p \int \prod_{j=1}^k \big|K_h(x_j-x_i)\big|
 \big(\mM[f](x_i)\big)^{p-k-1} f(x_i) \rd x_i
\\*[2mm]
& := g_k^{(1)}(x_1, \ldots, x_k)+
 g^{(2)}_k(x_1, \ldots, x_k).
\end{align*}
For the first term on the right hand side we obtain
\begin{align*} \rE_f |g_k^{(1)}(X_1, &\ldots, X_k)|^2 \leq c_1 \|K\|_\infty^{2(p-k)}
 \int \big(\mM[f](x_i)\big)^{2(p-k)}
 \\
 &\;\;\;\;\;\;\;\;\;\;\;\;\quad\times \;\bigg[\prod_{\substack{j=1\\ j\ne i}}^k \int K^2_h(x_j-x_i) f(x_j)\rd x_j\bigg]
f(x_i)\rd x_i
\\
&\leq \frac{c_1 \|K\|_\infty^{2p-2}}{V_h^{k-1}} \int \big(\mM[f](x)\big)^{2p-k-1} f(x)\rd x
\\
&\leq \frac{c_1 \|K\|_\infty^{2p-2}}{V_h^{k-1}} \int \big(\mM[f](x)\big)^{2p-k}\rd x \leq
\frac{c_2 \|K\|_\infty^{2p-2}}{V_h^{k-1}} \int [f(x)]^{2p-k} \rd x.
 \end{align*}
The expectation of the squared second term is bounded as follows:
\begin{align*}
 & \rE_f|g_k^{(2)}  (x_1, \ldots, x_k)|^2
 \\
 & \leq
 c_4 \|K\|_\infty^{2(p-k-1)}
 \iint
 \bigg\{
 \prod_{j=1}^k \int K_h(x_j-y)K_h(x_j-z) f(x_j)\rd x_j\bigg\}
 \\
 &\qquad\qquad\qquad\qquad\times \big(\mM[f](y)\big)^{p-k}\big(\mM[f](z)\big)^{p-k}\rd y\rd z
 \\*[2mm]
 &\leq \frac{c_3 \|K\|_\infty^{2p-2}}{V_h^{k-1}}
 \iint  {\bf 1}\{z-y\in I_{2h}(0)\} \big(\mM[f](y)\big)^{p}\big(\mM[f](z)\big)^{p-k}\rd z\rd y
 \\
 & \leq \frac{c_4 \|K\|_\infty^{2p-2}}{V_h^{k-1}}\int
\big(\mM[f](z)\big)^{p-k} \mM\big[\mM^p[f]\big](z)\rd z
\\
&\leq \frac{c_5 \|K\|_\infty^{2p-2}}{V_h^{k-1}}\int
 [f(z)]^{2p-k}\rd z,
\end{align*}
where we have used (\ref{eq:mu-M}).
Finally,
\begin{align*}
&\rE_f |g_p(X_1, \ldots, X_p)|^2\leq c_6 \int \bigg[\prod_{j=2}^p \int K_h^2(x_j-x) f(x_j)\rd x_j\bigg] f(x)  \rd x
\\
&\;\leq \frac{c_6\|K\|_\infty^{2p-2}}{V_h^{p-1}} \int  \big(\mM[f](x)\big)^{p-1} f(x) \rd x
\\
&\leq \frac{c_7\|K\|_\infty^{2p-2}}{V_h^{p-1}} \int \big(\mM[f](x)\big)^{p}\rd x\leq \frac{c_7\|K\|_\infty^{2p-2}}{V_h^{p-1}} \int [f(x)]^p\rd x.
\end{align*}
Thus, we obtain
\begin{align}
{\rm var}_f \big[\hat{T}_{2,h}\big] \leq  C_2 \|K\|_\infty^{2p-2}\sum_{k=1}^p
 \frac{1}{n^k V_h^{k-1}} \int [f(x)]^{2p-k}\rd x.
\label{eq:var-T2}
\end{align}
The assertion of the lemma follows now from (\ref{eq:var-T1}) and (\ref{eq:var-T2}).
%This completes the proof.
\epr

\paragraph{Proof of Lemma \ref{lem:inf-norm}}
 Let $K$ be a kernel satisfying Assumption \ref{ass2:kernel} with $\ell\geq \max_{j=1,\ldots,d}\beta_j$.
For any $\eta=(\eta_1, \ldots, \eta_d)$, $\eta_j>0$, $j=1,\ldots, d$ we have in view of Lemma \ref{lem:GL}
\begin{equation*}
%\label{eq:f-infty}
 \|f\|_{s}\leq \big\|B_\eta\big\|_{s} + \big\|S_\eta \big\|_{s}
 \leq  \sum_{j=1}^d \big\|B_\eta^{(j)}\big\|_{s} +
 \big\|K_\eta \ast f\big\|_{s}.
\end{equation*}
By the general form of Young's inequality with $1/q=1+1/s-1/p$
\begin{align*}
 \|K_\eta \ast f\|_{s} \leq \|f\|_p \|K_\eta\|_{q}
 \leq c_1\big(V_\eta\big)^{1/q-1}\|f\|_p=
 c_1\big(V_\eta\big)^{1/s-1/p}\|f\|_p.
\end{align*}
Furthermore, it follows from (\ref{eq:bias-2}) of
Lemma~\ref{lem:GL} with $s^*=s$ and $s_j=s$  that
\[
 \sum_{j=1}^d \big\|B^{(j)}_\eta\big\|_s \leq c_2\sum_{j=1}^d L_j \eta_j^{\gamma_j(s)},\quad \gamma_j(s)=\frac{\beta_j\tau(s)}{\tau(r_j)}.
\]
Therefore, for any $\eta_j>0, j=1,\ldots,d$
\[
 \|f\|_s\leq
  c_1\big(V_\eta\big)^{1/s-1/p}\|f\|_p +  c_2\sum_{j=1}^d L_j \eta_j^{\gamma_j(s)}.
\]
Putting $1/\gamma=\sum_{j=1}^d 1/\gamma_j(s)$ and
choosing
$\eta_1, \ldots, \eta_d$ from the equality
$$
\big(V_\eta\big)^{1/s-1/p}\|f\|_p=\sum_{j=1}^d L_j \eta_j^{\gamma_j(s)}
$$
 we come to the following
bound
\[
 \|f\|_s \leq c_3
 \bigg(\prod_{j=1}^d L_j^{\frac{1}{\gamma_j(s)}}\bigg)^
 {\frac{1/p-1/s}{1+\frac{1/p-1/s}{\gamma}}}
 \big\|f\big\|_p^{\frac{1}{1+\frac{1/p-1/s}{\gamma}}}.
 \]
Noting that
$$
1+\frac{1/p-1/s}{\gamma}=1+\bigg[\frac{1}{p}-\frac{1}{s}\bigg]\bigg[\frac{\tau(\infty)}{\beta\tau(s)}+ \frac{1}{\beta\omega\tau(s)}\bigg]=1+\frac{1/p-1/s}{\beta\tau(s)}
 =\frac{\tau(p)}{\tau(s)}
$$
we complete the proof.
\epr

\bibliographystyle{agsm}

\end{document}